\documentclass[11pt,a4paper]{amsart}
\usepackage{amsmath,amssymb,amsthm}
\usepackage[utf8]{inputenc}
\usepackage{graphicx}
\usepackage{color}

\swapnumbers

\theoremstyle{plain}
\newtheorem{thm}[equation]{Theorem}

\theoremstyle{definition}

\newtheorem{exmp}[equation]{Example}
\newtheorem{alg}[equation]{Algorithm}

\theoremstyle{remark}

\numberwithin{equation}{section}

\newcommand{\1}{\mathbf{1}}
\renewcommand{\d}{\mathrm{d}}
\newcommand{\E}{\mathbb{E}}

\renewcommand{\P}{\mathbb{P}}

\renewcommand{\phi}{\varphi}
\newcommand{\R}{\mathbb{R}}

\newcommand{\dist}{\mathrm{dist}}

\title[Simulation of Schr\"odinger equation]{Efficient simulation of Schr\"odinger equation with piecewise constant positive potential}

\author{Xuxin Yang}
\address{Xuxin Yang\\
Hunan First Normal University\\
Department of Mathematics\\
Changsha, Hunan 410205\\
P.R. of CHINA}

\author{Antti Rasila}
\address{Antti Rasila\\
Aalto University\\
School of Science\\
Department of Mathematics and Systems Analysis\\
P.O.Box 11100\\
FIN-00076 Aalto\\
FINLAND}

\author{Tommi Sottinen}
\address{Tommi Sottinen\\
University of Vaasa\\
Faculty of Technology \\
Department of Mathematics and Statistics\\
P.O.Box 700\\
FIN-65101 Vaasa\\
FINLAND}

\thanks{
This work is supported by the NNSF of China (No. 11571088), Hunan Provincial Natural Science Foundation of China
(14JJ7083), a Key Project Supported by Scientific Research Fund of Hunan Provincial Education Department (14A028), 
Scientific Research Fund of Hunan Provincial Education Department (14C0253), the Aid Program for Science and Technology Innovative
Research Team in Higher Educational Institutions of Hunan Province. A. Rasila was partially supported by Academy of Finland (No. 289576). T. Sottinen was partially funded by the Finnish Cultural Foundation (National Foundations' Professor Pool).  This work was done during T. Sottinen's research visit to Hunan First Normal University, Changsha, PRC.  T. Sottinen wishes to thank them for their hospitality.}

\keywords{
Brownian motion;
killing walk on spheres;
harmonic measure;
numerical algorithm;
Yukawa equation;
Schr\"odinger equation;
}

\subjclass[2010]{65C05; 68U20; 35Q40}

\begin{document}

\begin{abstract}
In this paper we introduce a new method for the simulation of a weak solution of the Schr\"odinger-type equation where the potential is piecewise constant and positive.  The method, called killing walk on spheres algorithm, combines the classical walk of spheres algorithm with killing that can be determined by using panharmonic measures.
\end{abstract}

\maketitle

\section{Introduction}

We investigate the simulation of the following special case of the time-independent Schrödinger equation
\begin{equation} \label{eq:schrodinger}
\frac12\Delta u(x) - q(x) u(x) = 0,
\end{equation}
where $D\subset \R^n$ , $n\ge 1$ is a bounded domain, $u \colon D \to \R$ is a piecewise two times continuously differentiable function and $q$ is of the form
\begin{equation}\label{eq:q-def}
q(x) = \sum_{m=1}^M \lambda_j \1_{D_m}(x), 
\end{equation}
where, for each $m=1,\ldots,M$, $\lambda_m$ is a non-negative constant, $D_m\subset D$ is a domain such that $D_m\cap D_\ell = \varnothing$ for $m\neq \ell$, and
$$
D= D\cap \big(\bigcup_{m=1}^M \overline{D_m}\big).
$$

This class of functions $q$ given by \eqref{eq:q-def} is motivated by the phenomenon called quantum tunneling in particle physics.  Also, the form \eqref{eq:q-def} can be motivated by being an approximation of general positive $q$.
 
Since we are considering the case where $q$ is discontinuous, we cannot obtain classical $u\in C^2(D)$ solution.  Instead, we consider weak solutions, i.e., $u$ is a solution of \eqref{eq:schrodinger} if
$$
\int_D u(x)\frac12\Delta\phi(x)\, \d x =
\int_D q(x)u(x)\phi(x)\, \d x 
$$
for all $\phi\in C_c^\infty(D)$.

\section{Panharmonic measures}

Let $W=(W_t)_{t\ge 0}$ be standard $n$-dimensional Brownian motion.  Let $\P^x$ be the probability measure such that, $\P^x$-almost surely, $W_0=x$.  Let $\E^x$ be the expectation associated with the measure $\P^x$. Let
$$
\tau_D = \inf\left\{ t>0\,;\, W_t \not\in D\right\},
$$
i.e., $\tau_D$ is the first exit time of the Brownian motion from the domain $D$. 
We assume that the domain $D$, and all the sub-domains $D_m$, are (Wiener) regular, i.e.,
$$
\P^y[\tau_{D}= 0] = 1 \quad\mbox{for all } y\in \partial D. 
$$
Then, by \cite[Theorem 4.7]{chung-zhao}, a bounded weak solution to the boundary value problem
\begin{equation}\label{eq:schrodinger-bvp}
\left\{
\begin{array}{rclll}
\frac12\Delta u(x) - q(x)u(x) &=& 0    &\mbox{on}& D, \\
                  u(x) &=& f(x) &\mbox{on}& \partial D,
\end{array}\right.
\end{equation}
where $f$ is a bounded Borel function, is given by
\begin{equation}\label{eq:kakutani-schrodinger}
u(x) = \E^x\left[ e^{-\int_0^{\tau_D} q(W_s)\, \d s} f\left(W_{\tau_D}\right)\right].
\end{equation}
Moreover, if $f$ is continuous, then \eqref{eq:kakutani-schrodinger} is the unique solution satisfying $u\in C(\bar D)$.

In the case where $q$ is of the form \eqref{eq:q-def} the formula \eqref{eq:kakutani-schrodinger} takes the form
\begin{equation}\label{eq:kakutani-q}
u(x) = \E^x\left[e^{-\sum_{m=1}^M \lambda_m T_m} f\left(W_{\tau_D}\right)\right],
\end{equation}
where $T_m$ is the total time the Brownian motion spends in the sub-domain $D_m$ before exiting the domain $D$ at time $\tau_D$.  Since the Brownian motion has the strong Markov property and the ``discounting term'' $e^{-\lambda_m T_m}$ can be considered as exponential killing, we can analyze the terms in the sum independently. This means that on the stochastic set $\{ t\ge 0; W_t\in D_m\}$ we can consider the corresponding Yukawa equation, or Brownian motion with exponential killing. This can be analyzed by using the panharmonic measures as in \cite{panharmonic}.
Now we recall some facts of panharmonic measures.

The \emph{harmonic kernel} is
\begin{equation}\label{eq:harmonic-kernel}
h^x(D;\d y, t) = \P^x\left[W_{\tau_D} \in \d y \,\mid\, \tau_D=t \right]
\frac{\d\P^x}{\d t}\left[\tau_D\le t\right].
\end{equation}
In \cite{panharmonic} it was shown that the harmonic kernel exists for all regular domains $D$ and that the harmonic measure can be written as
$$
H^x(D;\d y) = \int_{t=0}^\infty h^x(D;\d y, t)\, \d t.
$$
The \emph{$\lambda$-panharmonic measure} is
\begin{equation}\label{eq:panharmonic-measure}
H_\lambda^x(D; \d y) = \int_{t=0}^\infty e^{-\lambda t} h^x(D; \d y, t)\, \d t.
\end{equation}
It was shown in \cite{panharmonic} that all bounded solutions to boundary value problem to the Yukawa equation, i.e. \eqref{eq:schrodinger-bvp} with constant $q(x)\equiv \lambda>0$, can be represented via the $\lambda$-panharmonic measure as
\begin{equation}\label{eq:panharmonic}
u(x) = \int_{y\in\partial D} f(y)\, H_\lambda^x(D;\d y).
\end{equation}
It was also shown in \cite{panharmonic} that all the panharmonic measures $H_\lambda^x(D;\cdot)$, $\lambda\ge 0$ are equivalent.  Consequently, there exists Radon--Nikodym derivatives $Z_\lambda^x(D;\cdot)$ such that
$$
H_\lambda^x(D;\d y) = Z_\lambda^x(D;y)H^x(D;\d y).
$$
The harmonic measure corresponding to $\lambda=0$ is a probability measure. Indeed,
$$
H^x(D; \d y) =\P^x\left[W_{\tau_D}\in \d y\right].
$$
The panharmonic measure corresponding to case $\lambda>0$ is a sub-probability measure. Indeed,
\begin{equation}\label{eq:rn}
Z_\lambda^x(D;y) = \E^x\left[e^{-\lambda \tau_D} \,\Big|\, W_{\tau_D}=y\right]
\le 1.
\end{equation}
This Radon--Nikodym derivative and the classical harmonic measure play a central part in constructing the killing walk on spheres (KWOS) simulation algorithm in the next section.  The algorithm is an extension of the classical walk of spheres algorithm (WOS) due to Muller \cite{muller}.

\section{Construction of the algorithm}

The general formula \eqref{eq:kakutani-schrodinger} for the Schr\"odinger equation provides an obvious way to solve the boundary value problem \eqref{eq:schrodinger-bvp} by simulating Brownian particles.  Indeed, let $w^k_{x,\delta t}$, $k=1,\ldots,K$ be independent simulated Brownian trajectories starting from $x\in D$ with time-mesh $\delta t$. So, each $w^k_{x,\delta t}$ is an independent varying-length vector $(w^k_{x,\delta t}(0), w^k_{x,\delta t}(1), \ldots, w^k_{x,\delta t}(T^k_D))$, where $w^k_{x,\delta t}(0)=x$,
$$
w^k_{x,\delta t}(j) = w^k_{x,\delta t}(j-1) + \sqrt{\delta t}\,\xi_j^k,
$$
$\xi_j^k$'s are independent standard $n$-dimensional normal random variables, and $T^k_D$ is the first index $j$ such that $w^k_{x,\delta t}(j) \not\in D$. Set
\begin{equation}\label{eq:simu-trivial}
\hat u_{K,\delta t}(x) = \frac{1}{K} \sum_{k=1}^K \exp\left\{-\sum_{j=1}^{T^k} q\left(w^k_{x,\delta t}(j)\right)\, \delta t\right\}f\left(w^k_{x,\delta t}(T^k)\right).
\end{equation}
Then, by the strong law of large numbers and by the dominated convergence theorem, $u_{K,\delta t}(x) \to u(x)$ for all $x\in D$ with $\P^x$-probability 1 as $\delta t\to 0$ and $K\to\infty$.

The obvious problem in the general simulation method described above is that it is computationally very heavy.  However, since the solution of the Schr\"odinger equation \eqref{eq:kakutani-schrodinger} depends on the entire path of the Brownian motion, there seems to be no other way than take the $\delta t$ to be very small. Consequently, finding an efficient simulation algorithm for the general Schr\"odinger equation by using the Brownian motion seems a very challenging problem, indeed.

Let us then consider the special case of the Schr\"odinger equation, where $q$ is of the form \eqref{eq:q-def}.  This means that we have the Yukawa equation with different parameters in different regions. In this case, we can sometimes avoid simulating the Brownian particles in a fine time-mesh by killing the particle in the region with exponential clock.  

The idea of the simulation is comes from the following representation of the solutions, which follows from \cite[Corollary 2.11]{panharmonic}, the strong Markov property of the Brownian motion and the memoryless property of the exponential distribution.

\begin{thm}
Let $q$ be given by \eqref{eq:q-def}. Then the function $u$ in \eqref{eq:kakutani-schrodinger} can be give as
\begin{equation}\label{eq:killing-kakutani}
u(x) = \E^x\left[f(W_{\tau_D}); \tau^*>\tau_D\right],
\end{equation}
where $\tau^*$ is an independent exponential clock with intensity $\lambda_m$ on sub-domain $D_m$.
\end{thm}

Note that the formula \eqref{eq:killing-kakutani} depends on the entire path for the Brownian motion only through the killing time $\tau^*$.

Our estimator for $u(x)$ is
\begin{equation}\label{eq:simu-killing-kakutani}
\hat u_K (x) = \frac{1}{K}\sum_{k\in K^*(\lambda)} f(w_x^k(\tau_k)), 
\end{equation}
where $w_x^k$, $k=1,\ldots, K$ are independent simulations of the trajectories Brownian particles starting from point $x$, and the set $K^*(\lambda)\subset\{1,\ldots,K\}$ contains the particles that are not killed; $\tau_k$ is the termination-step of the algorithm.

The optimistic idea of how to individual trajectories $w_x=w_x^k$ are generated as follows: 
Set $w_x(0)=0$ and suppose $w_x(1),\ldots, w_x(j-1)$ are generated.
Suppose $w_x(j-1)\in D_m$, and not is near the boundary.  Suppose the $\lambda_m$-panharmonic Radon--Nikodym derivative  and the harmonic measure for the domain $D_m$ are known (this is the optimistic part).  Generate $w_x(j)$ on the boundary by using the harmonic measure $H_{\lambda_m}^{w_x(j-1)}(D; \d y)$.  Now $w_j(x)=y\in\partial D_m$. Kill the particle with probability $1-Z_{\lambda_m}^{w_x(j-1)}(D_m;y)$. If the particle dies, the algorithm terminates and the particle adds 0 to the sum \eqref{eq:simu-killing-kakutani}. If the particle survives and $y\in\partial D$, then the algorithm terminates and $f(w_x^k(j))$ is added to the sum \eqref{eq:simu-killing-kakutani}. If $y\in\partial D_{m'}\setminus \partial D$, then we have to simulate the particle with fine time-mesh in order to allow the particle to enter well inside the domain $D_{m'}$.  Once the particle is well inside a domain $D_{m'}$, we can generate a new particle by using the corresponding harmonic measure and the $\lambda_{m'}$-panharmonic Radon--Nikodym derivative as before.

The harmonic measures and panharmonic Radon--Nikodym derivatives are not known for arbitrary domain.  This limits the applicability of the idea above.  However,  for balls $B(x,r)\subset\R^n$ the Radon--Nikodym derivative $Z_\lambda^x(B(x,r);y)$, $\lambda>0$ is well-known. 
Indeed, by the rotational symmetry of the Brownian motion the Radon--Nikodym derivative is independent of $y$ and by the self-similarity of the Brownian motion 
$$
Z_\lambda^x(B(x,r); y) = \psi(\mu r),
$$
where $\mu = \sqrt{2\lambda}.$
The function $\psi$ is the Laplace transform of the first hitting time of a Bessel process on the level $1$ starting from $0$. This is well-known, see e.g. \cite[Theorem 2]{wendel} for $n\ge 2$:
\begin{equation}\label{eq:psi} 
\psi(\mu) = \frac{\mu^\nu}{2^\nu\Gamma(\nu+1)I_\nu(\mu)}, \quad \mu>0, 
\end{equation}
where 
$$
I_\nu(x) = \sum_{m=0}^\infty \frac{1}{m!\Gamma(m+\nu+1)}\left(\frac{x}{2}\right)^{2m+\nu}
$$
is the modified Bessel function of the first kind of order $\nu=n/2-1$.  The harmonic measure for balls is, due to the rotational symmetry of Brownian motion, simply the uniform distribution.  

This leads to the following simulation algorithm:

\begin{alg}[Killing walk on spheres (KWOS)]\label{alg:kwos}
\mbox{}
\begin{enumerate}
\item
Set $w_x(0)=x$ and suppose $w_x(1),\ldots, w_x(j-1)$ are generated.
\item
Generate $w_x(j)$ as follows:
\begin{enumerate}
\item
If $w_x(j-1)\in D_m$ and if $\dist(w_x(j-1),\partial D_m)<\varepsilon_2$ start generating the Brownian path with fine time-mesh $\delta t \ll \varepsilon_2$: generate
$$
w_x(j) = w_x(j-1) + \sqrt{\delta t}\xi(j),
$$  
where $\xi(j)$ an $n$-vector of independent standard normal random variables. Generate exponential killing: with probability $p(\delta t) = 1-e^{-\lambda\, \delta t}$. If killing occurs, the algorithm terminates and adds 0 to the sum \eqref{eq:simu-killing-kakutani}.
\item
If $w_x(j-1)\in D_m$ and if $\dist(w_x(j-1),\partial D_m)\ge\varepsilon_2$ generate $w_x(j)\in\partial B(w_x(j-1),r)$, $r=\dist(w_x(j-1), \partial D_m)$ uniformly.  Kill the particle with probability $1-\psi(\mu r)$, where $\psi(\mu r)$ is given by \eqref{eq:psi}. If killing occurs, the algorithm terminates and adds 0 to the sum \eqref{eq:simu-killing-kakutani}.
\item
If $\dist(w_x(j-1),\partial D)<\varepsilon_1$ let $w_x(j)$ be the projection of $w_x(j-1)$ to the boundary $\partial D$.  In this case the algorithm terminates and adds $f(w_x(j))$ to the sum \eqref{eq:simu-killing-kakutani}. 
\end{enumerate}
\end{enumerate}
\end{alg}

\section{Examples}

In this section, we present examples to motivate our algorithm and to illustrate its potential applications in one and two-dimensional settings. These examples were computed by using a simple implementation algorithms on Mathematica, and they mainly chosen from the point of view of visualisation. It should be noted that our approach is obviously more much attractive in higher dimensions, where many other methods for solving such are not available, or lead into excessive computation times\footnote{Time required for computing the examples in this section with a MacBook Air computer varied between seconds and minutes. Wolfram Mathematica 10.2 and very basic implementations of the algorithms with no performance optimisations were used in computations.}.

\subsection*{One-dimensional example}

Next we consider the equation \eqref{eq:schrodinger} in the case where $u\colon [a,b]\to \R$ and $[a,b]$ is a real interval, i.e., $\Delta u = u''$, and  $q(x)$ is as in \eqref{eq:q-def}. Then we obtain the impulsive differential equation
\begin{equation}
\label{onedim}
\frac{1}{2}u'' = \lambda(x)u,\quad \lambda(x)=\lambda_j\text{ for }x\in (x_{j-1},x_{j}),\quad j=2,\ldots,M,
\end{equation}
where $a=x_1<x_2<\cdots < x_{M-1}<x_M=b$, $M \ge 2$ is fixed, $u\in C^1([a,b])$, and $u''$ is piecewise continuous on $[a,b]$.

One should note that the theory of impulsive differential equations is much richer than the corresponding theory of differential equations
without impulse effects. For example, initial value problems of such equations may not, in general, possess any solutions at all even when the corresponding
differential equation is smooth enough, fundamental properties such as continuous dependence relative to initial data may be violated,
and qualitative properties like stability may need a suitable new interpretation. Moreover, a simple impulsive differential equation
may exhibit several new phenomena such as rhythmical beating, merging of solutions, and noncontinuability of solutions.    See e.g. \cite{Laksh} for more information about impulsive differential equations. However, these problems do not arise in the case of \eqref{onedim} as, the boundary value problem \eqref{eq:schrodinger-bvp} has a continuous weak solution, and it coincides piecewise with the classical solution.

The following algorithm is for the special case \eqref{onedim} with $a=x_1<x_2<x_3<x_4=b$ and
$$
q(x) = \lambda\1_{(x_2,x_3)}(x).
$$
To efficiently simulate the particles we combine the Gambler's ruin (i.e., the harmonic measure) outside the killing zone and the KWOS algorithm \ref{alg:kwos} inside the killing zone.  

\begin{alg}[Gambler's ruin with killing walk on spheres (GR-KWOS)]\label{alg:gamblers}
\mbox{}
\begin{enumerate}
\item Set $w(0)=x$ and suppose $w(1),\ldots, w(j-1)$ are generated.
\item Generate $w(j)$ as follows:
\begin{enumerate}
\item 
If $w(j-1)\in (x_1,x_2)$ then it will exit the domain in $x_1$ with probability $(x_2-w(j-1))/(x_2-x_1)$. The algorithm terminates and gives $u(x_1)$ to the sum \eqref{eq:simu-killing-kakutani}.  Otherwise $w(j)=x_2$ and we are entering the killing zone.
\item 
If $w(j-1)\in (x_3,x_4)$ then it will exit the domain in $x_4$ with probability $1-(x_4-w(j-1))/(x_4-x_3)$. The algorithm terminates and gives $u(x_4)$ to the sum \eqref{eq:simu-killing-kakutani}.  Otherwise $w(j)=x_3$ an we are entering the killing zone.
\item
If $w(j-1)\in [x_2,x_3]$ we are entering or in the killing zone. In this case use KWOS Algorithm \ref{alg:kwos}.
\end{enumerate}
\end{enumerate}
\end{alg}

The above algorithm can further be refined by using the harmonic measure and the $\lambda_m$-harmonic Radon--Nikodym derivative in all the domains $D_m$, $m=1,2,3$ (away from boundary).  Indeed, all the required formulas can be founded e.g. in \cite{borodin-salminen}.

Note that in the one-dimensional case, the function $\psi$ of \eqref{eq:psi} can be obtained from  \cite[(3.0.1)]{borodin-salminen}:
\[
\E^x[e^{-\lambda T}]  = \frac{\cosh( (b+a-2x)\sqrt{\lambda/2} ) }{ \cosh( (b-a)\sqrt{\lambda/2} )}, 
\]
where $T = \inf\{ t>0 ; W_t \not\in (a,b) \}$ and $W$ is the one-dimensional Brownian motion by using parameters $x=0$, $a=-r$, $b=r$, and $r>0$.

\begin{figure}[h]

\includegraphics[width=8cm]{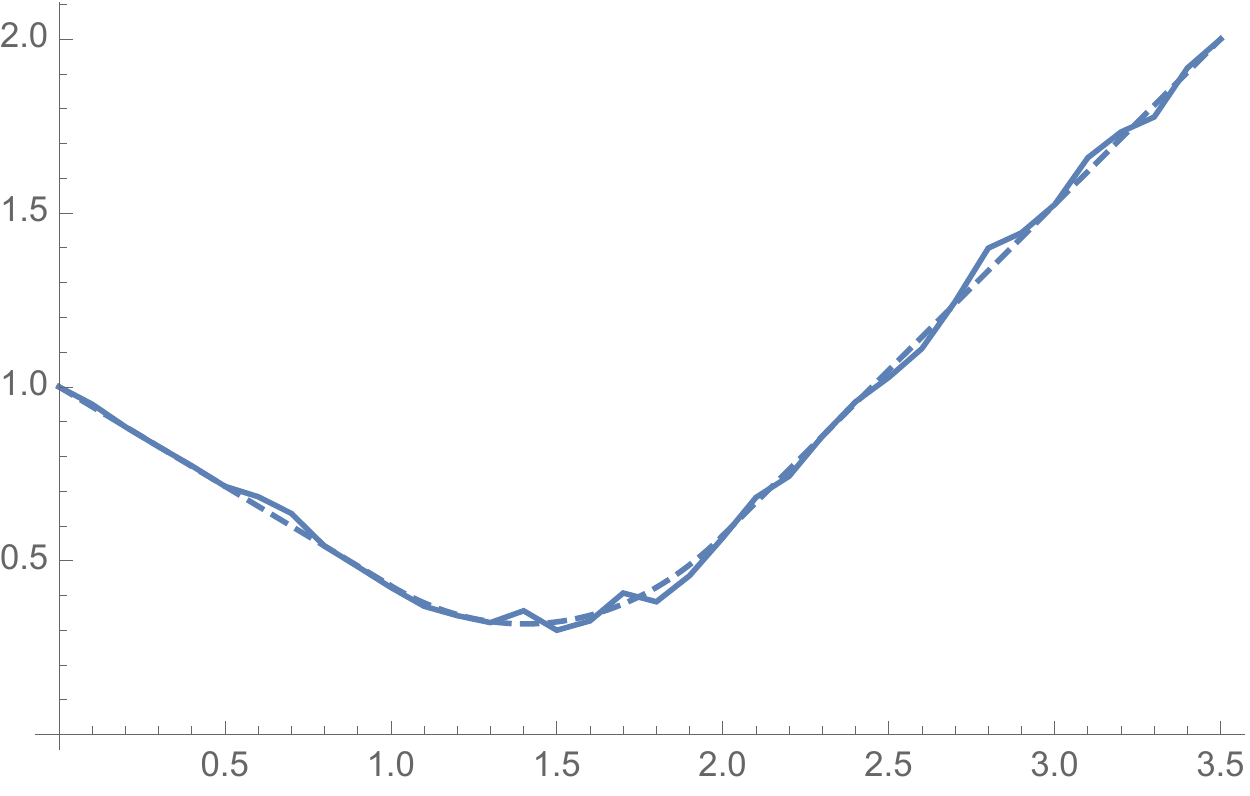}

\caption{The solution $u$ of the BVP of Example \ref{ex1} (dashed) and its approximation $\hat u_{1000}$ obtained by using GR-KWOS algorithm.}\label{ex1fig}
\end{figure}

\begin{exmp}
\label{ex1}
Let $x_1=0, x_2=1, x_3=2, x_4=3.5$, and consider the boundary value problem for the equation \eqref{onedim} with $u(x_1)=1, u(x_4)=2$, $\lambda_1 = \lambda_3=0$ and $\lambda_2=2$. Then the solution of the problem is given by
\[
u(x) \approx \left\{\begin{array}{rcl}
1-0.572084\, t & \text{for} & x\in (0,1),\\
0.00960025 e^{-2 t} (274.757 + e^{4t}) & \text{for} & x\in (1,2), \\
-1.33091 + 0.951688\, t & \text{for} & x\in (2,3.5).
\end{array}\right.
\]
The solution $u$, and approximate solution obtained through GR-KWOS Algorithm \ref{alg:gamblers}, are illustrated in Figure \ref{ex1fig}.
\end{exmp}

\subsection*{Two-dimensional examples}

Next we consider examples where the two dimensional potential function of the equation \eqref{eq:schrodinger} is computed by using Algorithm \ref{alg:kwos}. We start by computing the pure Yukawa case. Potential theory of the Yukawa equation has been studied by Duffin \cite{duffin1}. Stochastic methodology for studying Yukawa potentials was developed in \cite{panharmonic}.

\begin{figure}[h]

\includegraphics[width=8.2cm]{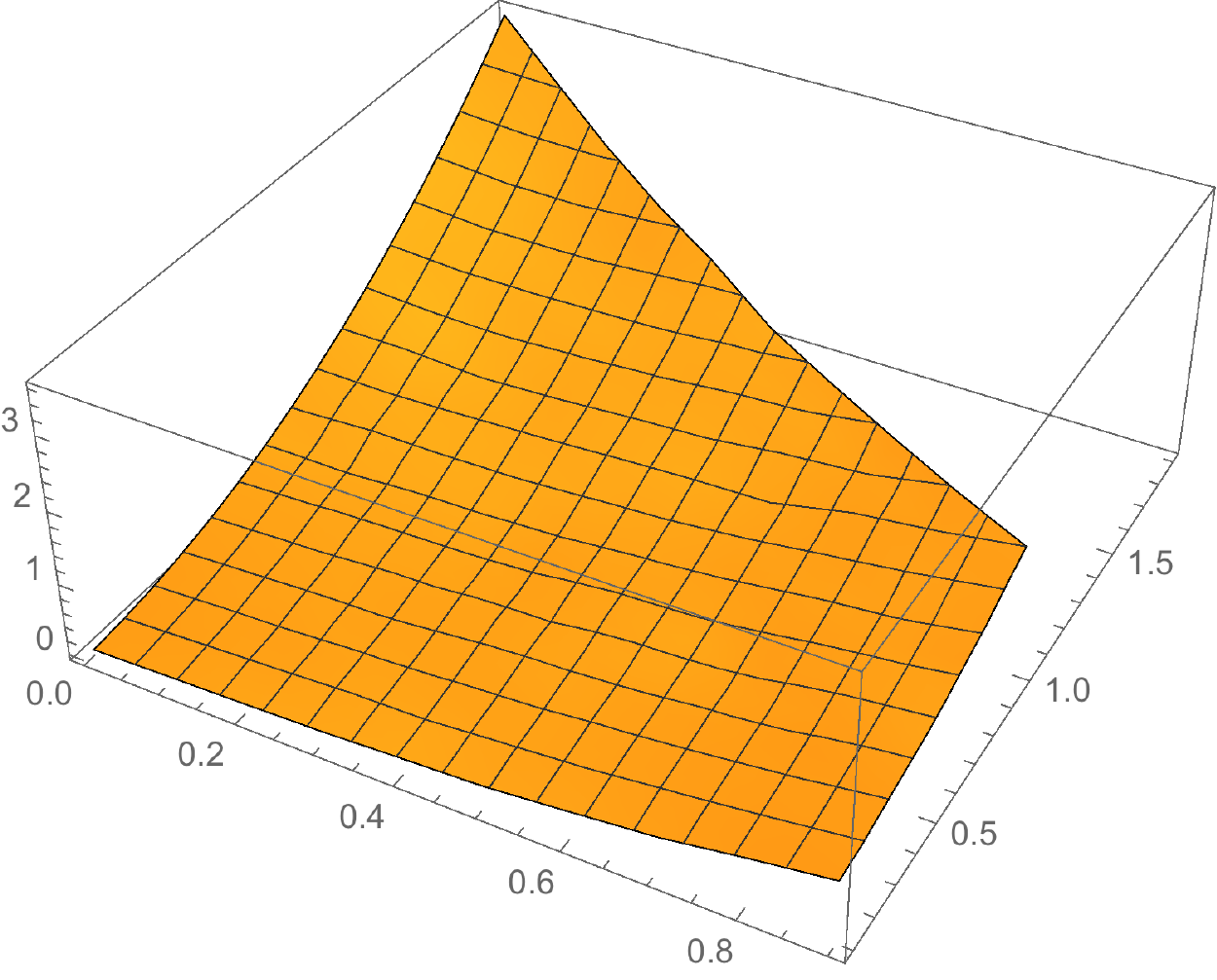}\quad
\includegraphics[width=4cm]{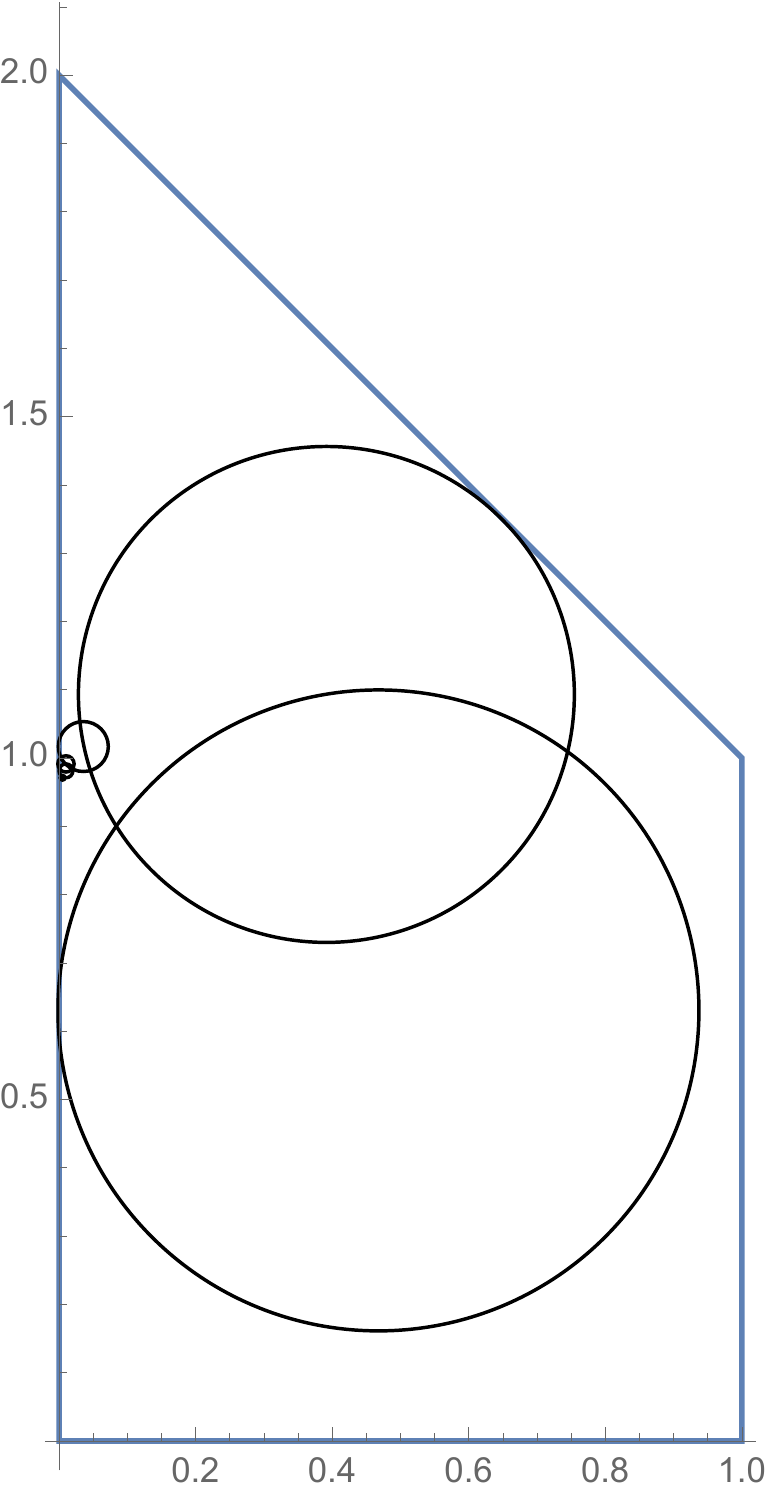}

\caption{An approximation of the solution $\hat u_{1000}$ of the BVP of Example \ref{ex2} with KWOS algorithm (left), and chain of disks used by KWOS algorithm in simulation of an individual path (right).}\label{ex2fig}
\end{figure}

\begin{exmp} (Yukawa equation)
\label{ex2}
Let 
$$
D=\big\{(x,y) ; 0<x<\min\{1,y+1\},\,0<y<2\big\}. 
$$
We solve the Yukawa equation, i.e., the equation
\eqref{eq:schrodinger} with constant $q\equiv 2$  with the boundary conditions given by $f(x,y)=x^3+y^2$. The approximate solution
$\hat u$ of \eqref{eq:simu-killing-kakutani}, where $K=1000$, and a chain of balls (disks) used in WOS style simulation of the path of one 
Brownian particle, are illustrated in Figure \ref{ex2fig}.
\end{exmp}

\begin{figure}[h]

\includegraphics[width=10cm]{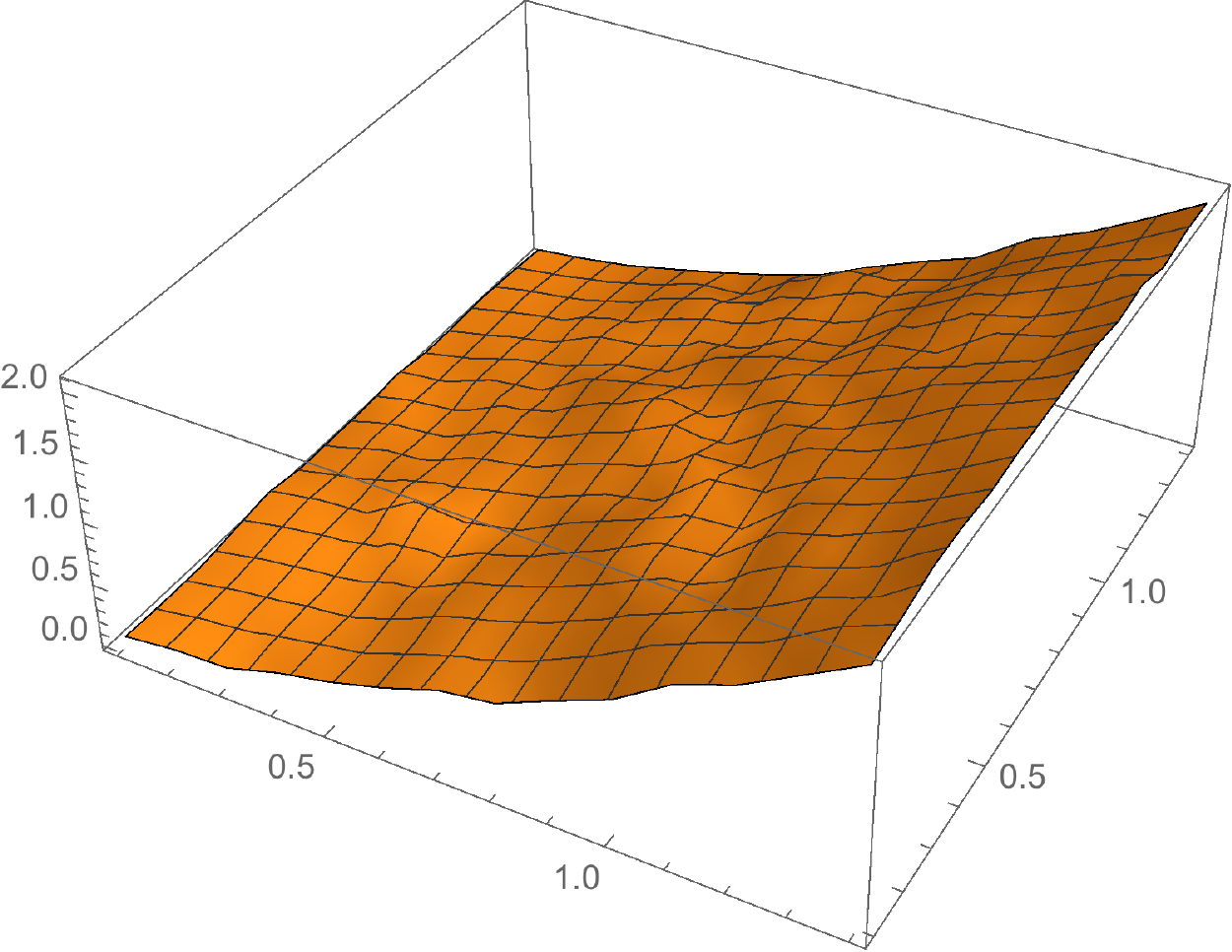}

\caption{The approximate  solution $\hat u_{250}$ of the BVP of Example \ref{ex3} with KWOS algorithm.}\label{ex3fig}
\end{figure}

\begin{exmp} (Mixed Laplace--Yukawa equation)
\label{ex3}
Consider the problem \eqref{eq:schrodinger} on the domain 
$$
D=\big\{ (x,y) ; 0< x < 1.5, \, 0<y<1.5\big\},
$$
where $q$ is as in \eqref{eq:q-def}, the boundary values are given by $f(x,y)=x^2$, 
\begin{eqnarray*}
D_1 &=& \big\{(x,y); 0<x<y,\,0<y<1\big\}, \\
D_2 &=& \big\{(x,y) ; 0< x < 1.5,\,x<y<1.5 \big\},
\end{eqnarray*}
$\lambda_1=0$ and $\lambda_2=2$.  The approximate solution
$\hat u_{250}$ of \eqref{eq:simu-killing-kakutani}, is illustrated in Figure \ref{ex3fig}. 
\end{exmp}

\section{Conclusion}

As a conclusion, the killing walk on spheres (KWOS) algorithm is a very simple tool to compute a solution to the Schr\"odinger equation with piecewise constant positive potential. The algorithm is based on the classical walk on spheres.  It avoids estimating the exit time from the ball by using the interpretation of killing.  If the potential has negative values, then the killing interpretation is no longer available and one needs to estimate the exit time by using, e.g., the recent walk on moving spheres (WOMS) algorithm due to Deaconu et al. \cite{woms,woms2}.

\end{document}